\documentclass[11pt,twoside]{amsart}

\usepackage{amssymb, palatino,graphicx, hyperref}
\usepackage[normalem]{ulem}
\usepackage[all]{xy}
\usepackage{xcolor}
\usepackage{ifthen}

\graphicspath{{./}{./Figures/}}

\newtheorem{dfn}{Definition}[section]
\newtheorem{thm}[dfn]{Theorem}
\newtheorem{rem}[dfn]{Remark}

\def\P{\ensuremath{\mathbb P}}

\def\R{\ensuremath{\mathbb R}}
\def\Z{\ensuremath{\mathbb Z}}

\def\C{\ensuremath{\mathbb C}}

\def\CP1{\ensuremath{\mathbb C \mathbb P^1}}

\def\Pn-1{\ensuremath{\P^{n-1}}}

\definecolor{red}{rgb}{.6,0,0}
\definecolor{green}{rgb}{0,.6,0}
\definecolor{darkgreen}{rgb}{0,0.3,0}
\definecolor{purple}{rgb}{0.5,0,0.5}
\definecolor{darkblue}{rgb}{0,0,0.7}
\definecolor{greenblue}{rgb}{0,0.4,0.5}
\definecolor{myblue}{HTML}{1685d6}
\definecolor{mypurple}{HTML}{b82e6c}

\newboolean{draft}

\newcommand{\cmt}[1]
{\ifthenelse {\boolean{draft}}
{{\sc \tiny \color{red} #1}}
{}}

\newcommand{\newbb}[1]
{\ifthenelse {\boolean{draft}}
{{\color{darkblue} #1}}
{#1}}

\newcommand{\newbbb}[1]
{\ifthenelse {\boolean{draft}}
{{\color{greenblue} #1}}
{#1}}

\newcommand{\nopost}[1]
{\ifthenelse {\boolean{draft}}
{{\color{cyan} #1}}
{}}

\newcommand{\maynopost}[1]
{\ifthenelse {\boolean{draft}}
{{\color{purple} #1}}
{}}

\newcommand{\margincmt}[1]
{\ifthenelse {\boolean{draft}}
{\marginpar{{\sc \tiny \color{red} #1}}}
{}}
\newcommand{\inred}[1]
{\ifthenelse{\boolean{draft}}{{\color{red} #1}}{#1}}

\newcommand{\new}[1]
{\ifthenelse {\boolean{draft}}
{{\color{green} #1}}
{#1}}

\newcommand{\neww}[1]
{\ifthenelse {\boolean{draft}}
{{\color{darkgreen} #1}}
{#1}}

\newcommand{\newb}[1]
{\ifthenelse {\boolean{draft}}
{{\color{blue} #1}}
{#1}}

\newcommand{\del}[1]
{\ifthenelse {\boolean{draft}}
{{\color{magenta} #1}}
{}}

\newboolean{details_on}

\newcommand{\details}[1]
{\ifthenelse {\boolean{details_on}}
{{\color{darkgreen} \tiny #1}}
{}}

\title{Nonrational polytopes and fans in toric geometry}
\author{Fiammetta Battaglia}
\author{Elisa Prato}
\thanks{The authors were partially supported by the PRIN Project ``Real and Complex Manifolds: Topology, Geometry and Holomorphic Dynamics" (MIUR, Italy) and by GNSAGA (INdAM, Italy).}

\newcommand{\bT}{\ensuremath{\mathcal{T} } } 

\newcommand{\calF}{\ensuremath{ \mathcal{F} } }

\setboolean{draft}{true}
\setboolean{details_on}{true}

\begin{document}
\maketitle
\vspace{0.6cm}
\begin{center}
\begin{minipage}[t]{11cm}
\small{
\noindent \textbf{Abstract.}
First, we examine the notion of nonrational convex polytope and nonrational fan in the context of toric geometry. 
We then discuss and interrelate some recent developments in the subject.
\medskip

\noindent \textbf{Keywords.}
toric variety, nonrational convex polytope, nonrational fan.
\medskip

\noindent \textbf{Mathematics~Subject~Classification:}
14M25, 52B20, 53D20.
}
\end{minipage}
\end{center}

\section*{Introduction}\label{introduzione} 

Toric varieties are a beautiful class of geometric objects, at the intersection of convex geometry and combinatorics on one side, and of algebraic and symplectic geometry on the other. 

The study and interest for toric varieties began in algebraic geometry with Demazure's foundational paper \cite{demazure}.
Some of the classical references on the subject are the article by Danilov \cite{danilov} and the books by
Fulton \cite{fulton} and Cox et al. \cite{cox}. From the symplectic perspective, the subject started with
Delzant's classification of symplectic toric manifolds \cite{delzant}, which is founded on the convexity 
theorem by Atiyah \cite{a} and Guillemin--Sternberg \cite{gs}.
Standard references for this viewpoint are the books by Audin \cite{audin}, Guillemin 
\cite{guilleminlibro} and Cannas de Silva \cite{cannas}.

The basic convex geometric objects that provide the starting point in classical toric geometry 
are {\em rational} convex polytopes and fans.

Our aim is to frame the notion of {\em nonrational} convex polytope and fan in the context of toric geometry. 
We give an historical account and then we describe, in the simplest possible way, how this notion has 
been recently interpreted by a number of authors who have dealt with the subject. The intent is to provide a 
unitary picture, a sort of dictionary, that makes it easier to move from one context to the other.

In Section~1, we recall the definitions of rational convex polytope and fan and view them in the toric geometric setting. 
We describe the fundamental starting convex data that are needed to extend toric geometry to the nonrational case. 
In Section~2, we describe a variant of the starting convex data. In Section~3, we illustrate the notions that were discussed in the previous sections with a number of examples. Finally, we dedicate Section~4 to the aforementioned dictionary.

\section{What is a nonrational convex polytope/fan: the fundamental triple}\label{triple}
A convex polytope $\Delta\subset (\R^n)^*$ is the convex hull of a finite number of points. Equivalently, 
it is the bounded intersection of finitely many closed half--spaces
$$\Delta=\bigcap_{j=1}^d\{\;\mu\in(\R^n)^*\;|\;\langle\mu,X_j\rangle\geq\lambda_j\;\},$$
where $X_1, \ldots X_d\in \R^n$, $\lambda_1,\ldots,\lambda_d\in \R$ and 
$d$ can be chosen to be exactly the number of codimension $1$ faces ({\it facets}) of $\Delta$ 
(see, for example, \cite[Theorem~1.1]{ziegler}). 
We assume, for simplicity, that $\Delta$ has maximal dimension $n$.
Remark that each vector $X_j$ is orthogonal to a facet of $\Delta$ and points towards its interior. 
We will be calling the vectors $X_1, \ldots X_d$ \textit{normals} for $\Delta$.
They are not unique, as each $X_j$, together with the corresponding $\lambda_j$, 
can be replaced by any positive scalar multiple. 

Convex polytopes are studied in combinatorics, but are also of fundamental importance in symplectic and algebraic geometry.
Think of the convexity theorem \cite{a,gs} and of geometric quantization in symplectic geometry. Or think of toric geometry in 
both algebraic and symplectic geometry. It is a crucial fact that the convex polytopes that appear in these classical geometric settings are all 
\textit{rational}. This means that they are always thought of together with a lattice. The precise definition of rational convex 
polytope goes as follows: a convex polytope $\Delta\subset (\R^n)^*$ is rational if there exists a lattice $L\subset\R^n$ such 
that the normals can be taken in $L$. Another crucial fact in toric geometry is that, for any rational convex polytope, 
there is a canonical choice of normals: each $X_j$ is taken to be the shortest possible vector in 
$L$, also known as {\em primitive} vector. Often, the primitive normals generate the lattice. We remark in passing that in 
symplectic toric geometry it is possible and interesting to consider also nonprimitive normals (see \cite{lertol}).

In the last few decades, there has been a growing interest in understanding how to make sense of toric varieties when the polytope is no longer rational. 

In 1999 \cite{pcras,p} the second author approached this question first by replacing the basic framework 
\textit{polytope--lattice--primitive normals}, which was clearly no longer suitable, with something more general. Let us explain.

\begin{figure}[h]
\begin{center}
\includegraphics[scale=0.5]{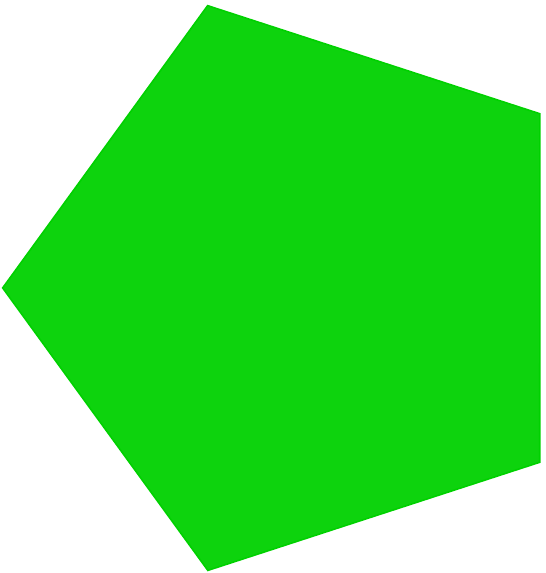} 
\caption{The regular pentagon}
\label{pentagon}
\end{center}
\end{figure}

\begin{figure}[h]
\includegraphics[scale=0.5]{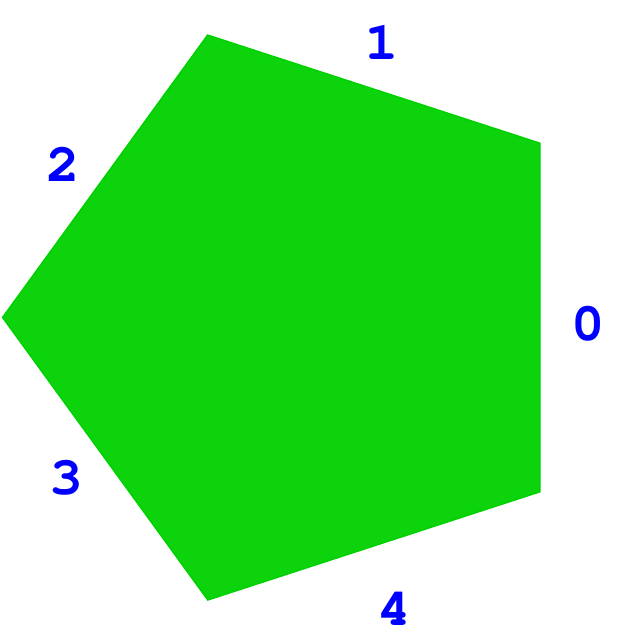} \quad\quad
\includegraphics[scale=0.3]{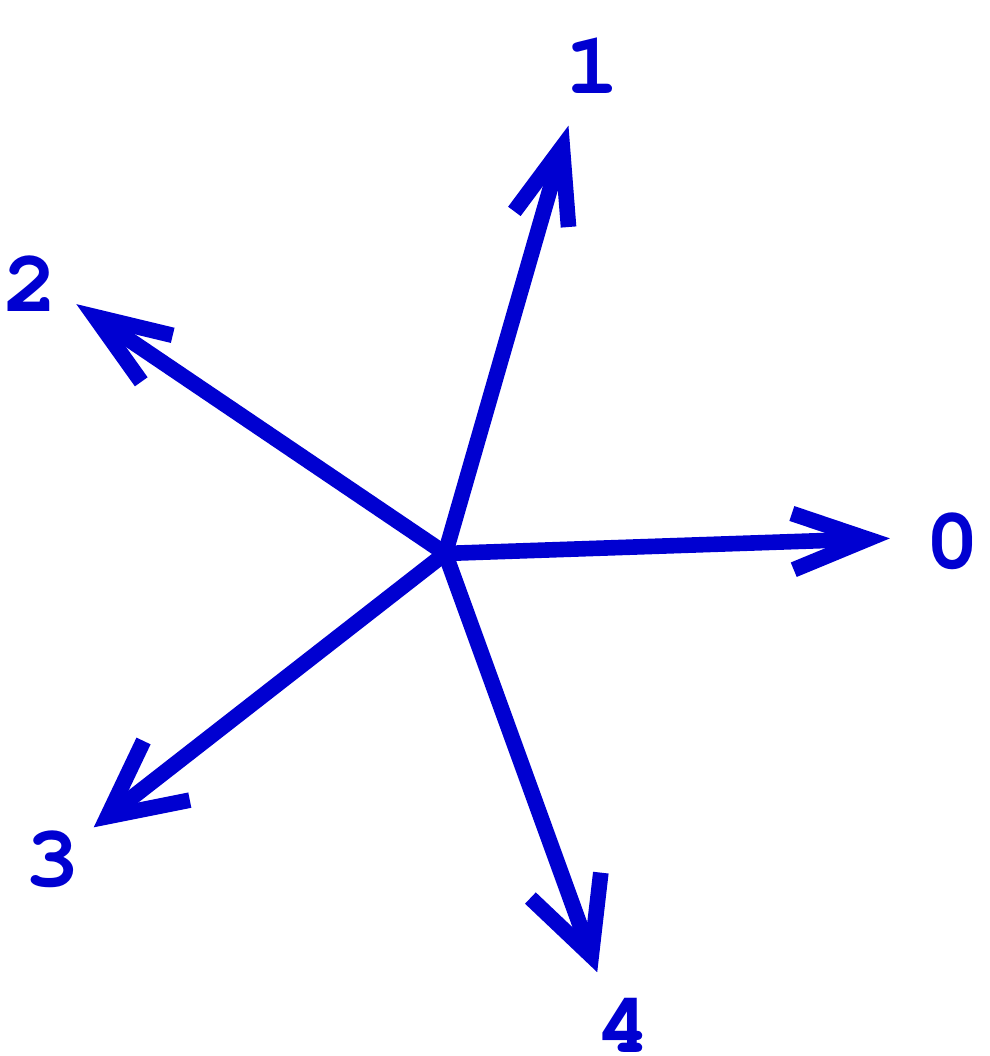}
\caption{The regular pentagon and the fifth roots of unity}
\label{pentagonandstar}
\end{figure}

The initial idea consisted in replacing the primitive vectors with any choice of normals, and the lattice with the 
$\Z$--span of these normals. The latter is a notion that was already well--known and of fundamental importance in the theory 
of quasicrystals and in the related theory of aperiodic tilings; it is called a \textit{quasilattice} \cite{mackay}. 
By definition, a quasilattice in $\R^n$ is the $\Z$--span of a set of $\R$--spanning vectors; 
notice that a quasilattice is a lattice if, and only if, these vectors form a basis of $\R^n$.
Take, for example, the regular pentagon (see Figure~\ref{pentagon}). 
It can be easily verified that it is not a rational polytope.
Consider the fifth roots of unity
$$
\begin{array}{l}
Y_0=(1,0)\\
Y_1=(\cos{\frac{2\pi}{5}},\sin{\frac{2\pi}{5}})\\
Y_2=(\cos{\frac{4\pi}{5}},\sin{\frac{4\pi}{5}})\\
Y_3=(\cos{\frac{6\pi}{5}},\sin{\frac{6\pi}{5}})\\
Y_4=(\cos{\frac{8\pi}{5}},\sin{\frac{8\pi}{5}}).
\end{array}
$$
Their opposites are normals for the pentagon (see Figure~\ref{pentagonandstar}).
The $\Z$--span of these vectors is a quasilattice, $Q_5$, that is dense in $\R^2$. We remark
that this particular quasilattice underlies the study of two aperiodic tilings of the plane discovered by Penrose:  
the rhombus tiling and the kite and dart tiling (see Section~\ref{examples}).

Now, two important remarks. First of all, the quasilattice generated by the chosen set of normals can be replaced by any other quasilattice that contains those normals. Take, for example, the polytope $[0,1]\subset\R^*$. 
It is rational with respect to the lattice $\Z\subset\R$, but also with respect to any lattice 
$a\Z$, $a\in\R$ (naturally isomorphic to $\Z$). However, we may also 
consider any quasilattice $Q=\Z+a\Z$, with $a$ irrational. 
Notice that this allows to consider rational convex polytopes in a nonrational setting. This turns out to be a natural choice in certain applications, some of which are described in Subsection~\ref{rationalinnonrational}. 
It also allows to perform symplectic cutting in an arbitrary direction \cite{cut}.

Secondly, one can consider any other set of normals, provided they are contained in the quasilattice, allowing even more freedom. 

This taken into account, it became convenient to define the following notion:
given a quasilattice $Q$, a convex polytope is said to be \textit{quasirational} with respect to $Q$ if the normals 
can be chosen in $Q$. Clearly, a convex polytope is quasirational with respect to a lattice 
if, and only if, it is rational. 
\begin{rem} Notice that, unlike rationality, quasirationality is not at all a restrictive requirement. In fact, a convex polytope is always quasirational with respect to the quasilattice that is generated by any set of normals. Think of the regular pentagon: it is not rational, however, it is quasirational with respect to the quasilattice $Q_5$. 
\end{rem}	
We are now ready to recall the notion of {\em fundamental triple}; it is the triple given by 
$$(\Delta,Q,\{X_1,\ldots,X_d\}),$$
where $\Delta\subset(\mathbb{R}^n)^*$ is any convex polytope, $Q\subset\mathbb{R}^n$ is any quasilattice with respect to which 
$\Delta$ is quasirational, and $\{X_1,\ldots,X_d\}$ is a choice of normals for $\Delta$ in $Q$. We remark that it is not required that $\{X_1,\ldots,X_d\}$ span the quasilattice.
The fundamental triple effectively replaces the \textit{polytope--lattice--primitive normals} triple of the rational case.
Once the triple is fixed, one can extend the classical geometric procedures for constructing toric varieties from polytopes.
For convex polytopes that are \textit{simple}, that is, when each vertex is the intersection of exactly $n$ facets, 
one gets a significant class of \textit{quasifolds}.

Quasifolds are highly singular spaces that are locally the quotient of a manifold modulo the action of a countable group. If the countable groups are all finite we get orbifolds, if they are all trivial we get manifolds. As it happens for manifolds, even for quasifolds the local models are required to be mutually compatible and thus form an atlas.
Quasifolds are naturally endowed with the usual geometric objects such as vector fields, differental forms and in particular symplectic structures. For the formal definition of quasifold and related notions, in the real and complex setting, we refer the reader to \cite{pcras,p,kite,cx}. 

Going back to nonrational toric geometry, 
the real and complex tori $\mathbb{R}^n/L$ and $\mathbb{C}^n/L$ of the rational case are naturally replaced by the quotients $\mathbb{R}^n/Q$ and $\mathbb{C}^n/Q$, which are abelian groups and quasifolds that are referred to as \textit{quasitori} \cite{pcras,p,cx}. 
We are now ready to recall the following basic results. Let $\Delta$ be a simple convex polytope. From \cite[Theorem~3.3]{p} we have
\begin{thm}\label{thmsp}
For each fundamental triple $(\Delta,Q,\{X_1,\ldots,X_d\})$, there exists a compact, connected $2n$--dimensional quasifold $M$, endowed with a symplectic structure and an effective Hamiltonian action of the quasitorus $\R^n/Q$ such that, if $\Phi\,\colon M\rightarrow(\R^n)^*$ is the corresponding moment mapping, then $\Phi(M)=\Delta$. 
\end{thm}
Moreover, from \cite[Theorem~2.2]{cx}
\begin{thm}\label{thmcx}
For each fundamental triple $(\Delta,Q,\{X_1,\ldots,X_d\})$, there exists a compact, connected, $n$--dimensional complex quasifold $X$, endowed with a holomorphic action of the complex quasitorus $\C^n/Q$ having a dense open orbit. 
\end{thm}
Finally, from \cite[Theorem~3.2]{cx}
\begin{thm}\label{thmkhl}
For each fundamental triple $(\Delta,Q,\{X_1,\ldots,X_d\})$, the symplectic quasifold $M$ of Theorem~\ref{thmsp} is equivariantly diffeomorphic to the complex quasifold $X$ of Theorem~\ref{thmcx}. The symplectic and complex structures are compatible and hence define a K\"ahler structure on $M\simeq X$.
\end{thm}
The space $M\simeq X$ is known as the \textit{toric quasifold} associated with the triple $(\Delta,Q,\{X_1,\ldots,X_d\})$. As in the rational case, it is explicitly constructed by means of symplectic and complex quotients.
In the proofs of Theorems~\ref{thmcx} and \ref{thmkhl} one sees that it is endowed with two beautiful finite atlases that generalize the standard complex affine and symplectic toric atlases of the rational case; their charts are modeled on $\C^n\simeq\R^{2n}$ modulo the action of countable subgroups of the standard torus $\R^n/\Z^n$.
 
The case of general convex polytopes was addressed by the first author, who showed that the resulting space $M$, which is naturally even more singular, is stratified by toric quasifolds \cite{stratif-re,stratif-cx}.

We remark that, unlike what happens in the rational setting, usually we do not have canonical choices for quasilattices and normals. 
Sometimes, however, the general geometric setup suggests natural choices. Some instances are described in Section~\ref{examples}. This is also the case when performing symplectic reduction and symplectic cutting in the nonrational toric setting \cite{cut,riduzione}.

The notion of rational convex polytope can be naturally expressed in terms of the rationality of its normal fan. The same applies also to quasirationality. Let us first recall the definition of fan, which is the central convex object in the theory of toric varieties in algebraic geometry.
A \textit{fan} $\Sigma$ in $\R^n$ is a collection of cones such that each nonempty face of a cone in $\Sigma$ is itself a cone in $\Sigma$ and such that the intersection of any two cones in $\Sigma$ is a face of each \cite{ziegler}.
The one--dimensional cones are said to be the \textit{generating rays} of the fan. The \textit{normal fan} $\Sigma_{\Delta}$ of the 
convex polytope $\Delta$ is the fan whose generating rays are inward pointing and orthogonal to the polytope facets and such that 
there is an inclusion--reversing bijection between cones in $\Sigma_{\Delta}$ and faces of $\Delta$. It is a {\em complete} fan, namely the union of its cones is $\R^n$. 
\begin{figure}[h]
	\includegraphics[scale=0.3]{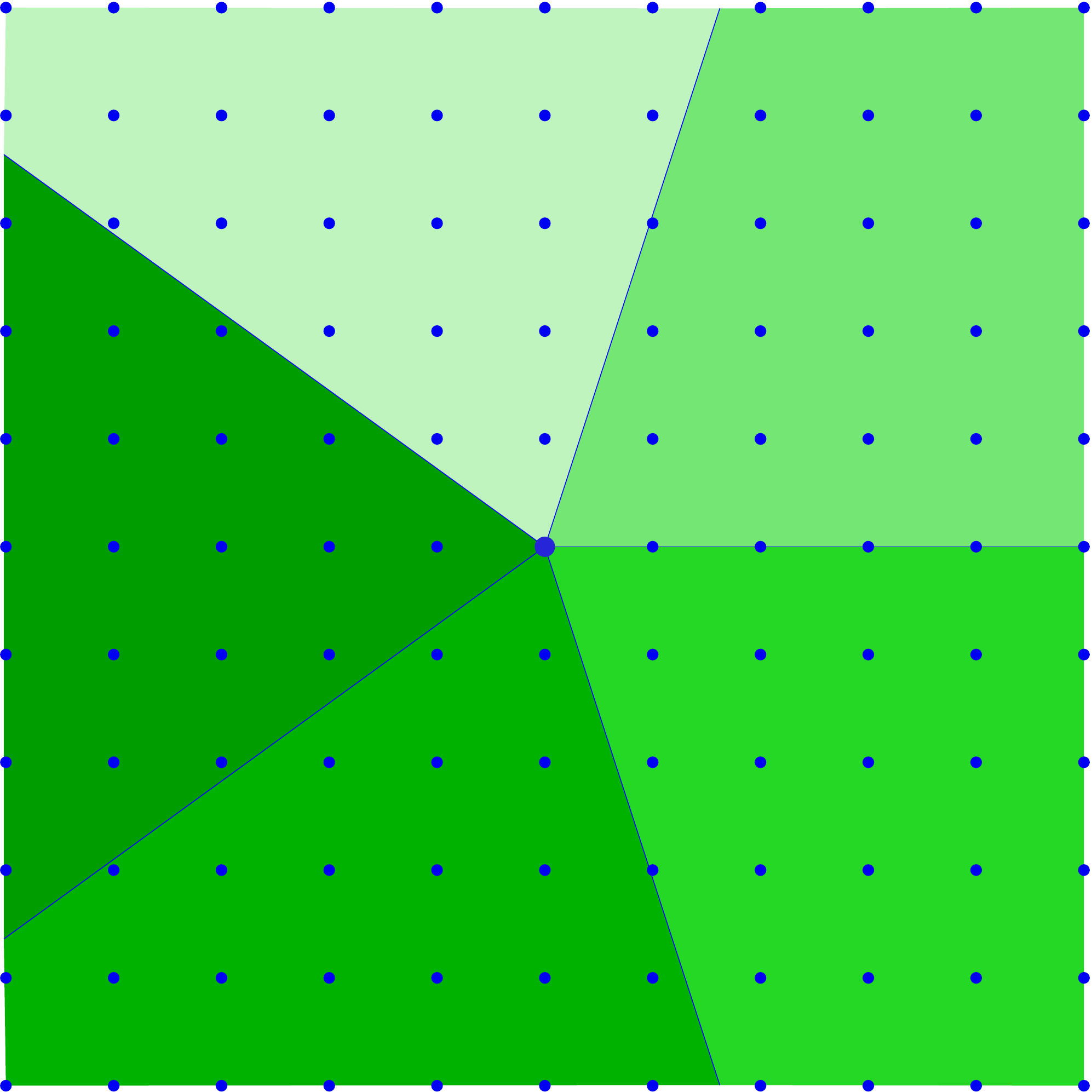}
	\caption{The normal fan of the regular pentagon}
	\label{pentagonfan}
\end{figure}
Moreover, a fan is said to be {\em simplicial} 
if each of its cones is simplicial, namely spanned by linearly independent vectors. 
Notice that a convex polytope is simple if, and only if, its normal fan is simplicial.
The rationality/quasirationality of the polytope corresponds to the rationality/quasirationality of its normal fan. 
We recall, in fact, that a fan in $\R^n$ is said to be
rational if there exists a lattice $L\subset\R^n$ which has non--empty intersection with each generating ray.
Similarly, we say that a fan in $\R^n$ is quasirational with respect to a quasilattice $Q\subset\R^n$ 
if $Q$ has non--empty intersection with each generating ray.

For any fan $\Sigma$, we can still introduce the triple
$$(\Sigma,Q,\{X_1,\ldots,X_d\}),$$ where $Q\subset\mathbb{R}^n$ is any quasilattice with respect to which 
$\Sigma$ is quasirational and where the vectors $\{X_1,\ldots,X_d\}$ are generators of the fan rays in the quasilattice $Q$. When the fan is complete and simplicial the complex construction of Theorem~\ref{thmcx} applies verbatim.
When, in addition, the fan is {\em polytopal,} namely when $\Sigma$ is the normal fan $\Sigma_{\Delta}$ of a convex polytope $\Delta$, 
the choice of such a $\Delta$
endows $M$ with the symplectic, and hence K\"ahler, structures of Theorems~\ref{thmsp} and \ref{thmkhl}.
We can draw the following commutative diagram
$$\xymatrix{(\Delta,Q,\{X_1,\ldots,X_d\})\ar[rd]\ar[d]&\\(\Sigma_{\Delta},Q,\{X_1,\ldots,X_d\})\ar[r]&M.}$$

\section{The fundamental triple encoded in a triangulated vector configuration: the augmented triple}
\label{triangolazione}
As we have seen, toric quasifolds are constructed explicitly and share key features with their rational counterparts. 
However, they are highly singular topological spaces. 
It is therefore natural to ask if there exists a framework that allows to work with smooth objects. 
With this motivation in mind, the first author, jointly with Zaffran \cite{bz1}, 
introduced in 2011 the idea of viewing toric quasifolds as leaf spaces of compact, complex, holomorphically foliated manifolds. 
This development was built on, and inspired by, two previous articles: the above--mentioned article by the second author
on nonrational toric geometry \cite{p}, and the article \cite{MV} by Meerssemann--Verjovsky.
In the latter, simplicial projective toric varieties were already
viewed, in the classical rational setting, as leaf spaces of LVM manifolds, a large class of compact, complex, 
non--K\"ahler manifolds \cite{ldm,M}, admitting a holomorphic foliation \cite{LN,M}.
The viewpoint developed in \cite{bz1} naturally brought a new perspective on the convex geometric data as well: the fundamental 
triple was encoded in a triangulated vector configuration, a well known and studied convex object \cite{DL-R-S}. 
Let us recall from \cite[Section~2]{bz1} what a triangulated vector configuration is, how a triple is encoded there and, finally, why this convex datum 
is instrumental in the construction of LVMB manifolds -- a generalization of LVM manifolds \cite{bosio}.

An {\em odd, balanced, triangulated vector configuration} is given by a pair $(V,\bT)$, where $V=(X_1,\ldots,X_p)$ is an 
ordered list of vectors in $\R^{n}$, allowing repetitions, that is {\em balanced}, namely 
$\sum_{i=1}^{p} X_i=0$, and {\em odd}, namely $p-n=2m+1$. A subset $\tau$ of
$\{1,\ldots,p\}$ is a simplex when the vectors indexed by $\tau$ are linearly independent. 
The cone generated by these vectors is called $\text{cone}(\tau)$.
A {\em triangulation} $\bT$ of a configuration $V$ is a collection of simplices
satisfying the following conditions:
\begin{enumerate}
\item If $\tau\in\bT$ and $\tau'\subset \tau$ then $\tau'\in\bT$
\item For all $\tau,\tau' \in\bT$, $\text{cone}(\tau) \cap \text{cone}(\tau')=\text{cone}(\tau \cap \tau')$
\item $\cup_{\tau\in\bT}\ \text{cone}(\tau) \supset \text{cone}(V)$.
\end{enumerate}
Remark that $(V,\bT)$ encodes
\begin{itemize}
\item a simplicial fan, not necessarily polytopal: the union of the cones indexed by $\bT$
\item ray generators (the vectors indexed by $\bT$)
\item a quasilattice $Q=\text{Span}_\Z(X_1,\dots,X_p)$
\item a number of {\em ghost} vectors (those which are not indexed by $\bT$).
\end{itemize}
Viceversa, let $\Sigma$ be a simplicial fan. We can construct a triangulated vector configuration $(V,\bT)$ that encodes a given triple $(\Sigma,Q,\{X_1,\ldots,X_d\})$ as follows:
if $\text{Span}_\Z\{X_1,\ldots,X_d\}=Q$ and the vector configuration $(X_1,\ldots,X_d)$ is odd and balanced, we 
keep it as it is, otherwise we can add ghost vectors, so as to have a set of generators of the quasilattice $Q$ and a balanced, odd configuration. 
Notice that there are infinitely many choices of ghost vectors that comply with 
these conditions. We remark that, in the rational case, this procedure of adding vectors so as to have a set of generators of the 
lattice and a balanced, odd configuration, 
is already used in \cite{MV}, where additional vectors correspond to {\em indispensable points}.  The fan is complete if, and only if, the vectors in $V$ span $\R^n$.
In conclusion, the vector configuration $V$ takes care of the quasilattice and the vectors in the triple, while $\bT$ is determined by the fan combinatorics.

In short, we consider the {\em augmented triple} $$(V,\bT)=\text{triple}+\text{ghost vectors},$$ where we have chosen a set of generators for the quasilattice $Q$ that includes a set of ray generators. 
The technical conditions, for the vector configuration to be balanced and odd are not at all restrictive. They allow to obtain, from a pair $(V,\bT)$, the exact convex datum that produces an LVMB manifold. More precisely, Gale duality applied to $V$ gives
a configuration $\Lambda=(\Lambda_1,\ldots,\Lambda_p)$ of points in affine space $\C^m$
(viceversa, Gale duality applied to $\Lambda$ determines $V$ up to automorphisms  \cite[Section~1.2]{rome}). On the other hand, the combinatorial datum $\bT$ yields a combinatorial datum $\bT^*$, which is a {\em virtual chamber} of $\Lambda$.
In turn, each pair $(\Lambda,\bT^*)$ 
determines a compact, complex, holomorphically foliated manifold $(N,\calF)$, of complex dimension $n+m$, where $m$ is the dimension of the leaves \cite[Section~2.2.4]{bz1}. 
We can draw the following diagram:
\begin{equation}
\label{diagramma}
\xymatrix{(V,\bT)\ar@{-->}[r]\ar[d]&(\Lambda,\bT^*)\ar[r]&(N,\calF)\ar[d]\ar[d]\\
(\Sigma,Q,\{X_1,\ldots,X_d\})\ar@{-->}@<1ex>[u]\ar[rr]&&M}
\end{equation}
Two arrows are dashed as they are not maps, in the sense that the target object is not uniquely determined. This implies that, to a given fundamental triple there corresponds a whole family of LVMB manifolds of different dimensions. When the fan is polytopal, $N$ is an LVM manifold and the leaf space of each member of this family is exactly the corresponding complex toric quasifold $M$ (see \cite[Section~2.3.1]{bz1} and \cite[Theorem~2.1]{rome}). 
Polytopality of the fan $\Sigma$ can be expressed in terms of the above--mentioned convex objects (for more details see \cite{bz1} and references therein). In terms of the corresponding LVM manifold, polytopality implies that the foliation is transversely K\"ahler \cite{LN,M}. The converse is also true; the proof, by Ishida \cite{I1}, is based on a convexity theorem in this context. 
In the polytopal case, diagram (\ref{diagramma}) can be constructed in the symplectic setting: $N$ turns out to be a presymplectic manifold, while the symplectic quasifold $M$ is given by $N$ modulo the action of a connected abelian group \cite{rome}. Taking the augmented triple ensures that the presymplectic manifold $N$ is even dimensional and that the group is connected.

Subsequently, other important convexity theorems were proved, for transversely symplectic foliations \cite{ratiu,LS} and for \'etale symplectic stacks \cite{HS}.

\section{Examples}\label{examples}
As far as examples go, two situations arise naturally. In the first, we have examples of convex polytopes/fans that are 
nonrational. In the second, we have rational convex polytopes/fans inside of a geometric context where it is interesting, and 
sometimes downright necessary, to replace the lattice by a suitable quasilattice.
\subsection{Purely nonrational}
\begin{figure}[h]
\includegraphics[scale=0.3]{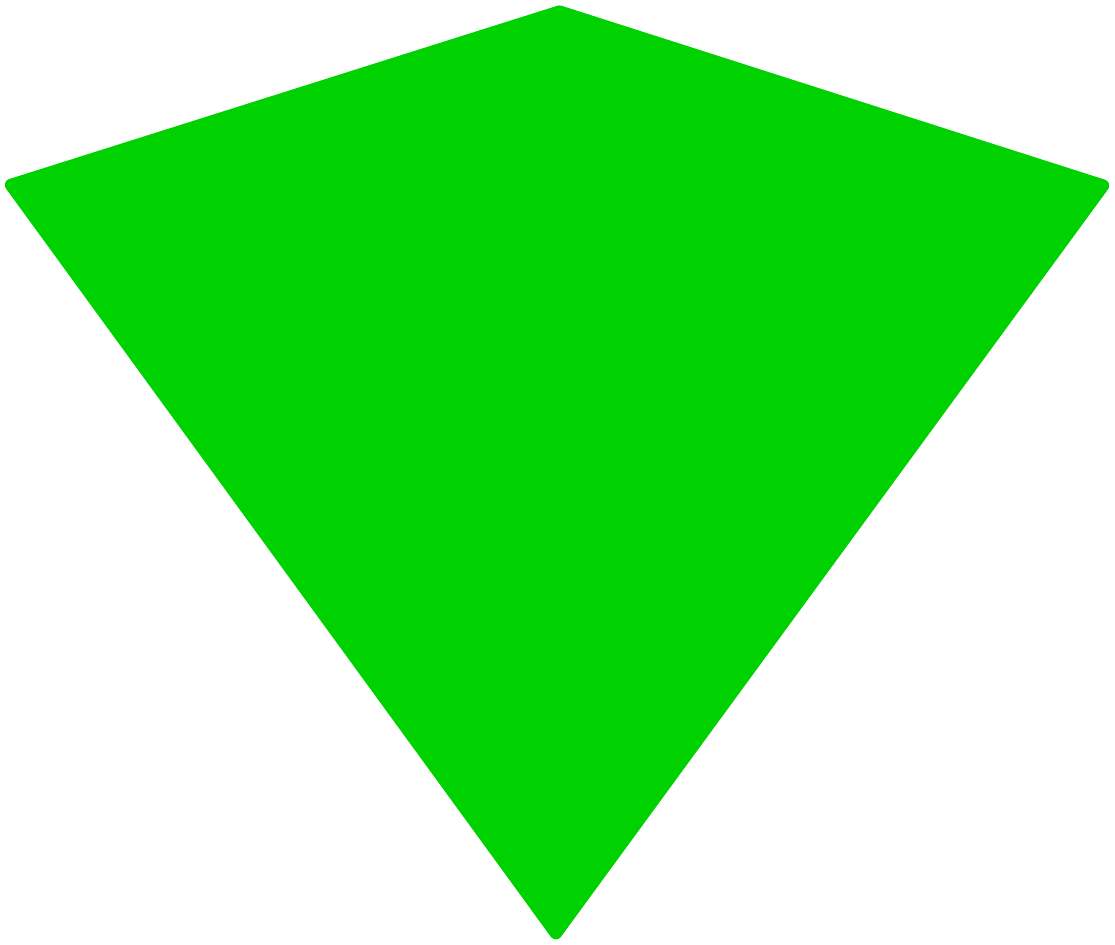}
\caption{The Penrose kite}
\label{kite}
\end{figure}

The first example that comes to mind of a purely nonrational convex polytope is the regular pentagon that we discussed in Section~\ref{triple}. 
An even simpler example is given by the Penrose kite, which, together with the dart, 
has been used by Penrose to construct aperiodic tilings of the plane \cite{pentaplexity}. 
The kite is the quadrilateral pictured in Figure~\ref{kite}. 
Three of its angles equal $\frac{2\pi}{5}$, while the other equals $\frac{4\pi}{5}$. Moreover,
its long edge is $\phi$ times its short edge, where
$\phi=\frac{1+\sqrt{5}}{2}$ is the \textit{golden ratio}.
The regular pentagon and the Penrose kite are actually closely related; in fact, the kite can be obtained from the regular pentagon 
via a standard construction (see, for example, \cite{kite}). Moreover, the kite is quasirational with respect to the same 
quasilattice $Q_5$ that we introduced for the pentagon.
\begin{figure}[h]
\includegraphics[scale=0.3]{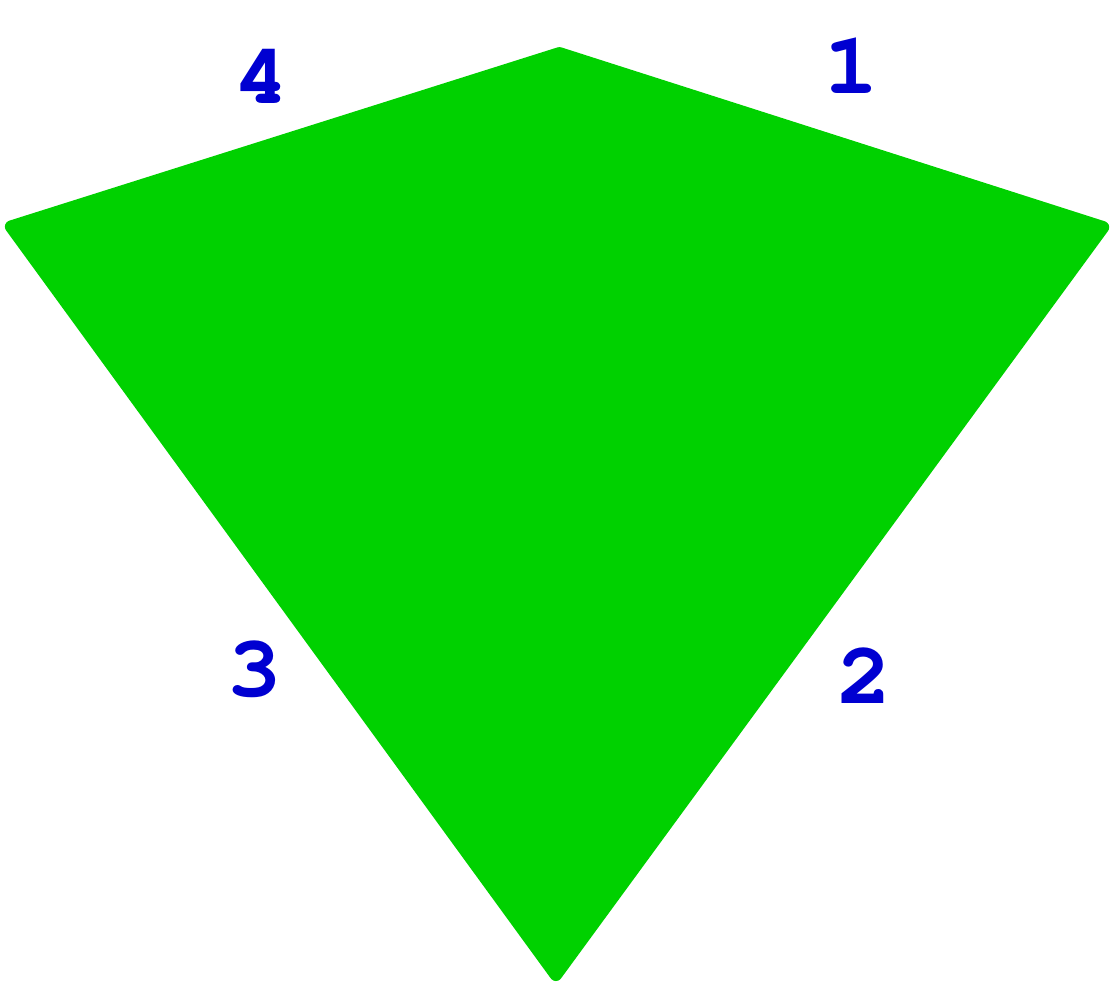} \quad\quad
\includegraphics[scale=0.3]{star.pdf}
\caption{Quasirationality of the Penrose kite}
\label{kiteandstar}
\end{figure}
In fact, a natural choice of normals for the kite is given by
$-Y_1$, $Y_2$, $-Y_3$, $Y_4$ (see Figure~\ref{kiteandstar}). 
These vectors generate the four rays of the corresponding normal fan. 
A vector configuration here is given by $V=(-Y_1,Y_2,-Y_3,Y_4,Y_0)$;
it is odd and balanced, with triangulation, in term of maximal simplices, given by
$\bT=\{\{1,4\},\{4,3\},\{3,2\},\{2,1\}\}$.

More details on the Penrose kite
from the symplectic toric viewpoint can be found in \cite{kite}.

Other interesting and elementary examples of nonrational convex polytopes are given by the regular 
dodecahedron and the regular icosahedron, only the first of which is simple. We refer the reader 
to \cite{dodecahedron,platonics} for a suitable choice of quasilattice and normals and 
for a description of the corresponding toric spaces. 

All of the above are closely related to the physics of quasicrystals, see for example \cite{ammann,senechal}.

We conclude by remarking that any simple nonrational convex polytope can be perturbed into a combinatorially equivalent rational one. This is not necessarily the case for nonsimple convex polytopes. See
\cite{grunbaum,ziegler_int} for the first counterexample, due to Perles. 

\subsection{Rational convex polytope in a nonrational setting}
\label{rationalinnonrational}
We have already seen in Section~\ref{triple} how the unit interval $[0,1]$, obviously rational, can be viewed as quasirational. The corresponding toric quasifold is a \textit{quasisphere} \cite{pcras, p,p4}. 

We illustrate two other examples where this is natural.
	
\begin{figure}[h]
\includegraphics[scale=0.3]{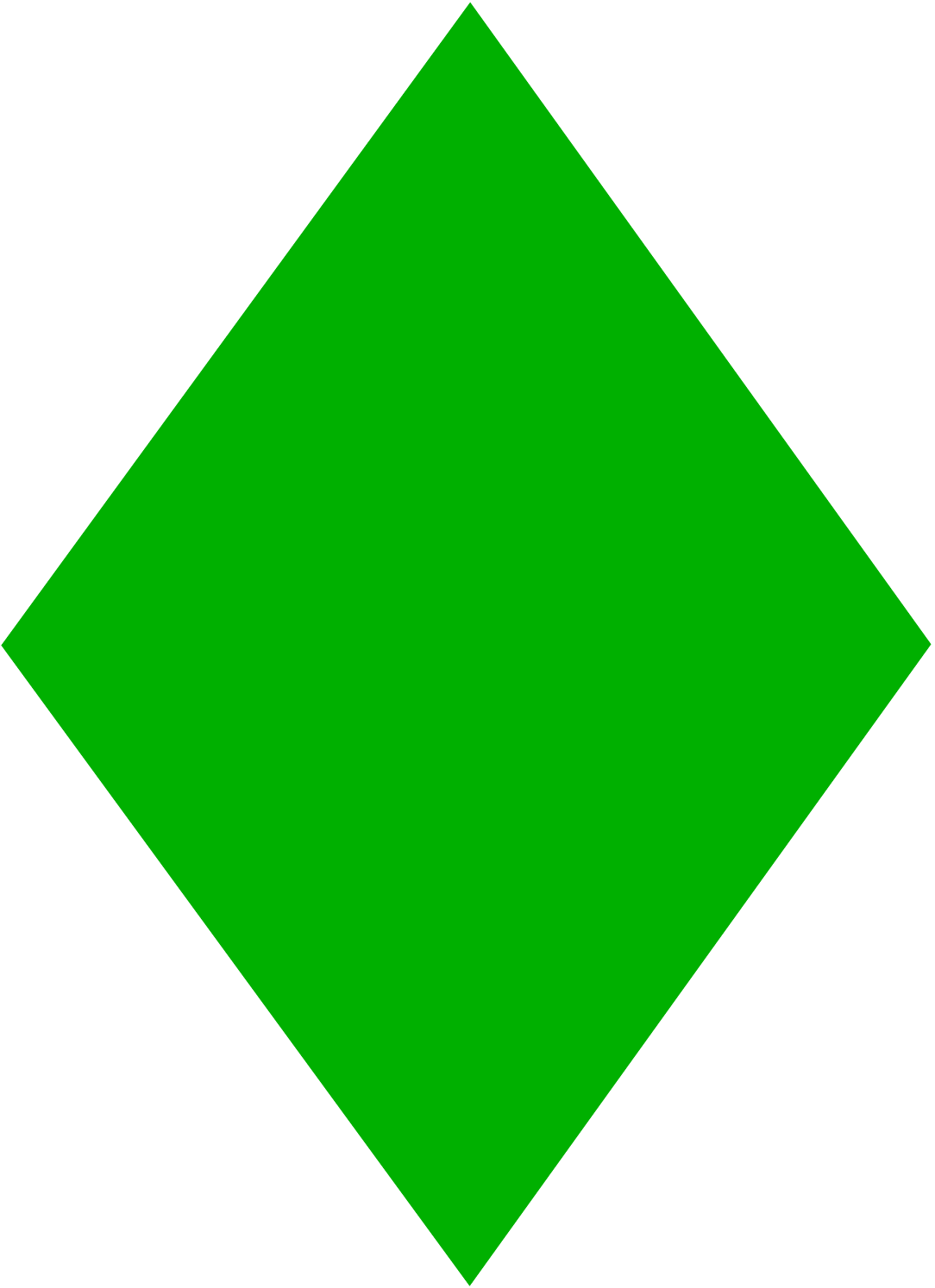}
\includegraphics[scale=0.3]{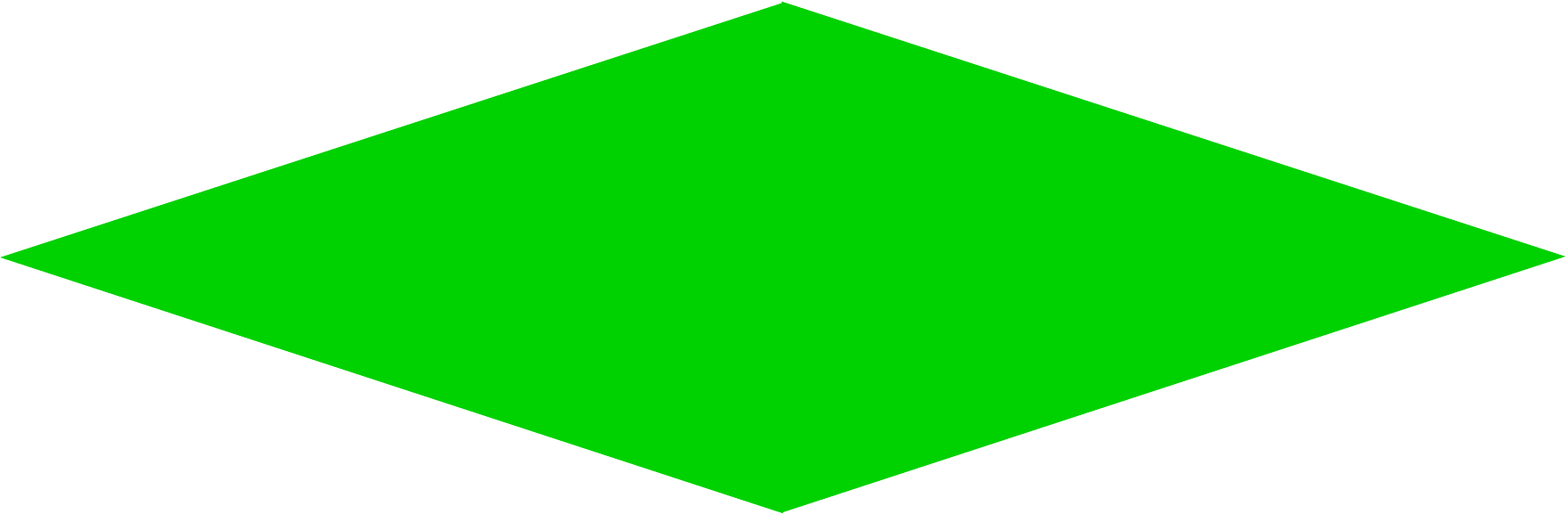}
\caption{The thick and thin Penrose rhombuses}
\label{rhombuses}
\end{figure}

\begin{figure}[h]
\includegraphics[scale=0.3]{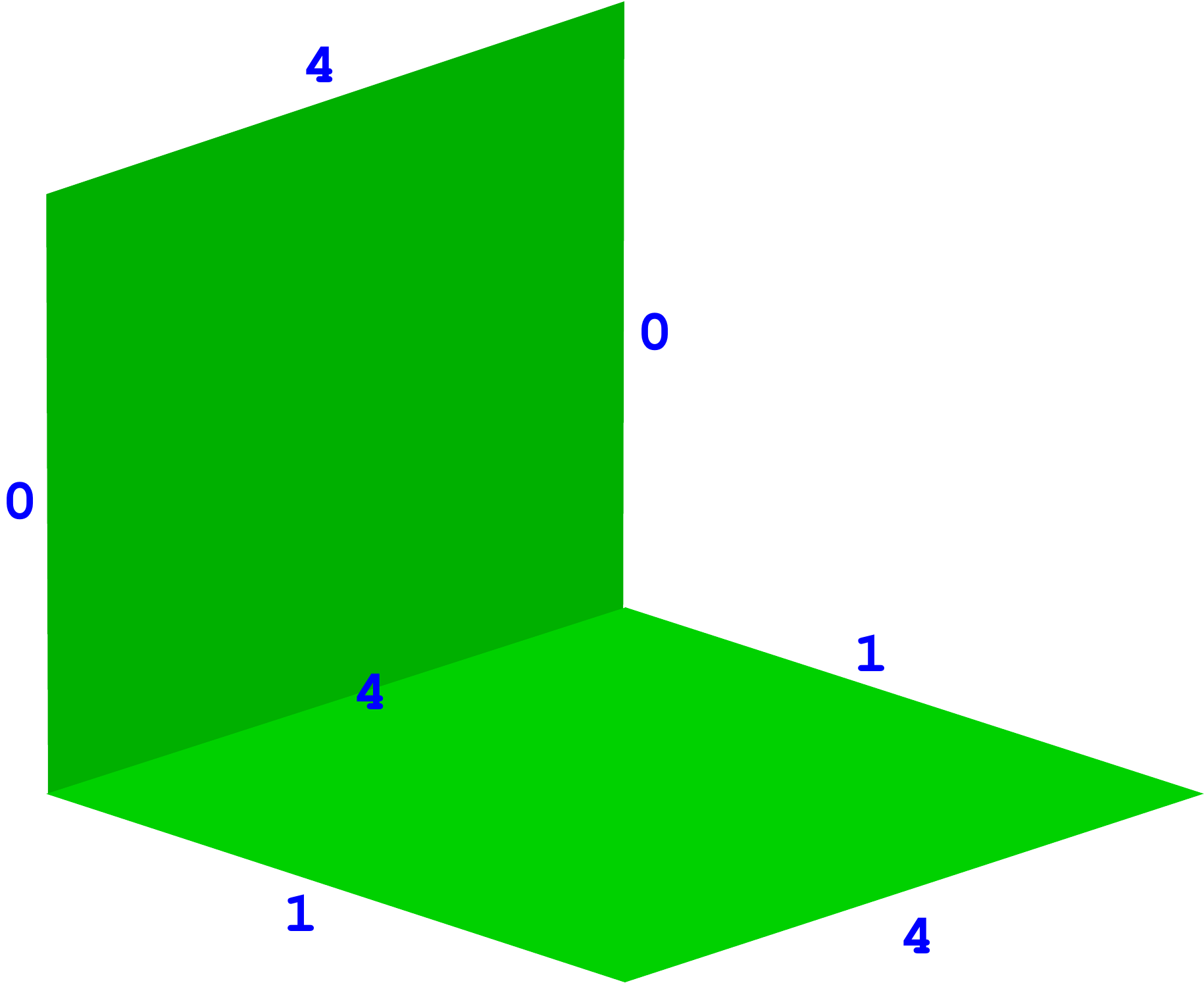} \quad\quad
\includegraphics[scale=0.3]{star.pdf}
\caption{Quasirationality of the Penrose rhombuses}
\label{rhombusesandstar}
\end{figure}

Let us consider first the Penrose rhombus tiling \cite{pentaplexity}. It is another fundamental
aperiodic tiling, whose tiles are given by two types of rhombuses, 
known as \textit{thick} and \textit{thin}. The thick rhombus has angles equal to $\frac{2\pi}{5}$ 
and $\frac{3\pi}{5}$, and the long diagonal is given by $\phi$ times the edge;
the thin rhombus has angles equal to $\frac{\pi}{5}$ 
and $\frac{4\pi}{5}$, and the edge is given by $\phi$ times the short diagonal
(see Figure~\ref{rhombuses}).
Each of them viewed individually is actually a rational convex polytope, but it is natural to want to consider a 
geometric setup that takes into account the entire tiling. In order to do so, they need to be viewed as quasirational with 
respect to the same quasilattice $Q_5$ that we have considered for the pentagon and the kite. 
In Figure~\ref{rhombusesandstar} we see normals $\pm Y_0$, $\pm Y_4$ for the thick rhombus
and $\pm Y_1$, $\pm Y_ 4$ for the thin one.
 To obtain vector configurations here we need to consider more vectors than we did in the case of the kite. 
For example, for the thick rhombus, we can take $V=(Y_0,Y_4,-Y_0,-Y_4,Y_1,Y_2,Y_3+Y_4+Y_0)$;
it is odd and balanced, with same triangulation as the kite. 
We refer the reader to \cite{rhombus} for further details on rhombus tilings from the symplectic viewpoint.

\begin{figure}[h]
\includegraphics[scale=0.5]{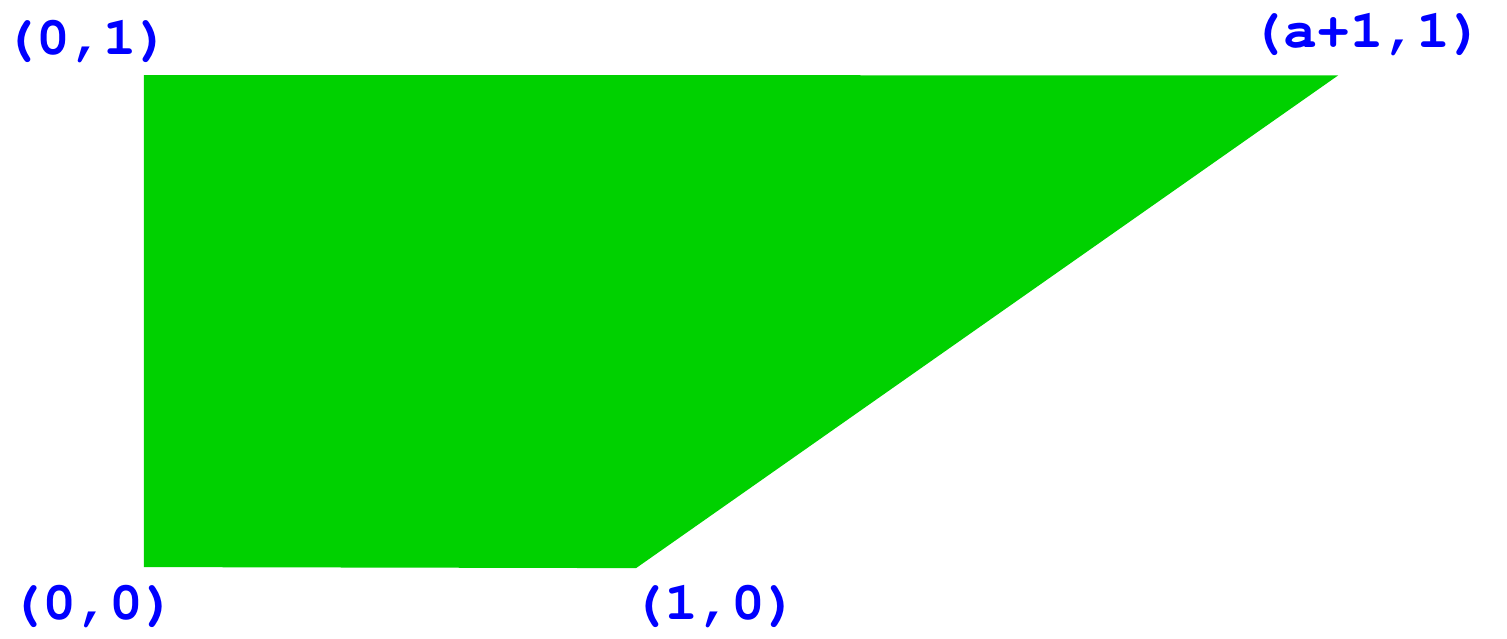}
\caption{The trapezoid $T_a$}
\label{quasihirze}
\end{figure}

\begin{figure}[h]
	\includegraphics[scale=0.5]{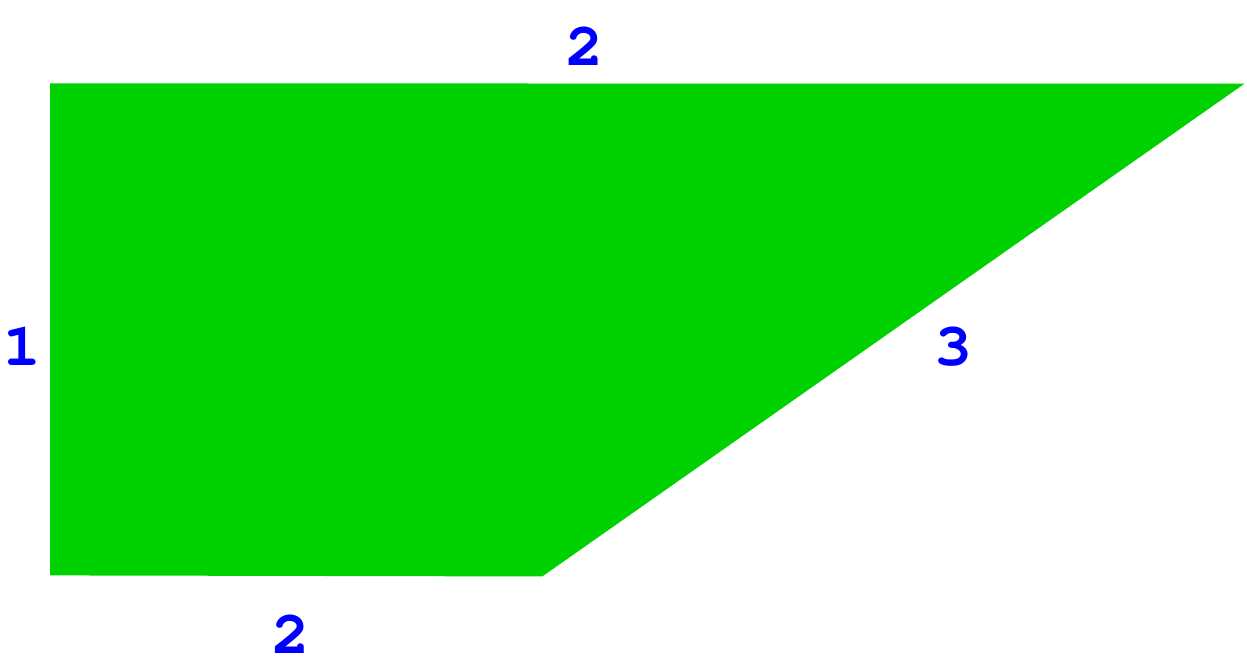} \quad\quad
	\includegraphics[scale=0.3]{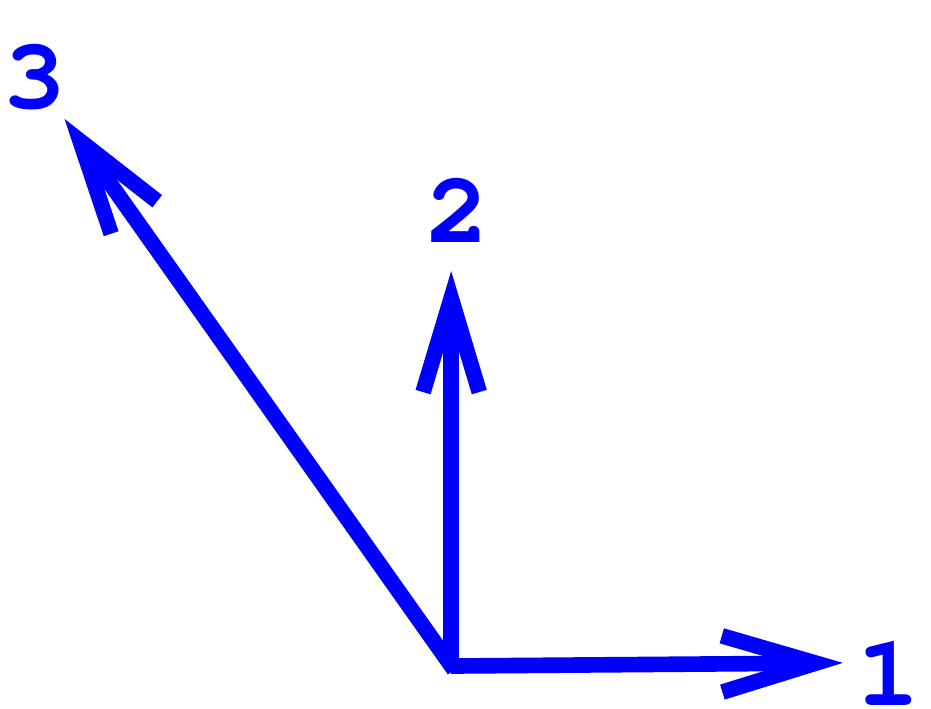}
	\caption{Quasirationality of the trapezoid $T_a$}
	\label{quasihirzeandfan}
\end{figure}

Other interesting examples arise when generalizing Hirzebruch surfaces to the nonrational setting.
Namely, consider, for any positive real number $a$, the trapezoid $T_a$ of vertices $(0,0)$, $(1,0)$, $(0,1)$ and $(a+1,1)$ 
(see Figure~\ref{quasihirze}).
When $a$ equals a positive integer $n$, we get the trapezoid $T_n$ that corresponds, in standard toric
geometry, to the Hirzebruch surface $\mathbb{H}^n$. We recall that this toric variety is constructed relatively 
to the standard lattice $\Z^2$ and to the primitive normals $(1,0)$, $\pm (0,1)$, $(-1,n)$.
For $a$ irrational, the trapezoid $T_a$, though rational with respect to the lattice that is generated by 
$(1,0)$ and $(0,a)$, is not rational with respect to the standard lattice $\Z^2$. 
If we want to consider a setup that yields, as a special case, 
the standard one for Hirzebruch surfaces, it is necessary to consider 
normals $(1,0)$, $\pm (0,1)$, $(-1,a)$ (see Figure~\ref{quasihirzeandfan}),
which span the quasilattice $Q_a=\Z\times (\Z+a\Z)\supseteq \Z^2$. Notice that, for $a$ rational, $Q_a$ is a lattice
and that, for $a=n$, this lattice equals $\Z^2$, as required. 
\begin{figure}[h]
	\includegraphics[scale=0.3]{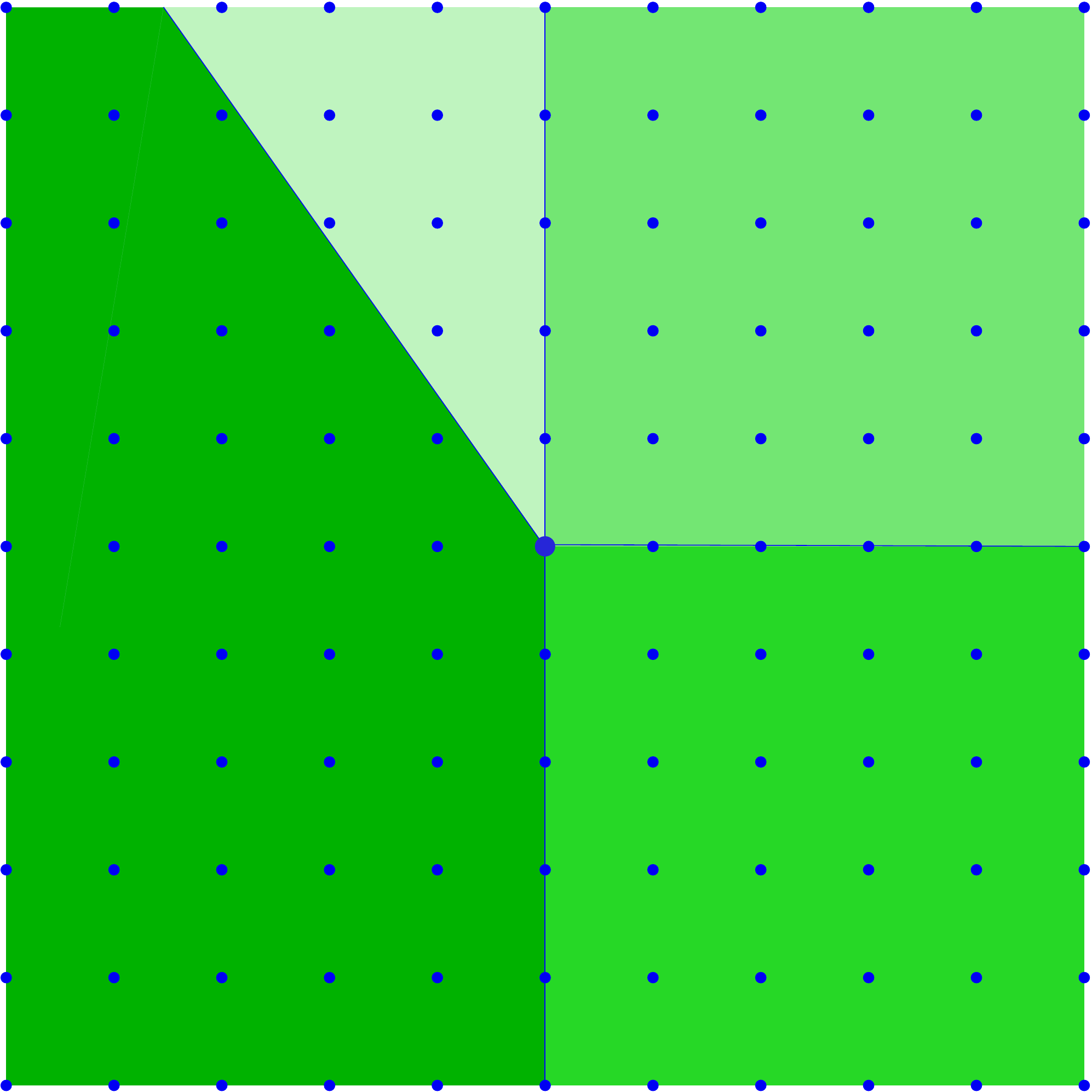}
	\caption{The normal fan of the trapezoid $T_a$}
	\label{quasihirzefantwo}
\end{figure}
The normal fan here is the complete fan in $\R^2$ whose generating rays are spanned by the four normals 
$(1,0)$, $\pm (0,1)$ and $(-1,a)$ (see Figure \ref{quasihirzefantwo}). 
A corresponding triangulated vector configuration is given by 
$(V_a,\bT)$, with 
$$V_a=((1, 0), (0, 1), (0,-1), (-1, a), (0,-a)),$$
$$\bT=\{\{1,2\},\{2,4\},\{3,4\},\{1,3\},\{1\},\{2\},\{3\},\{4\},\varnothing\}.$$ We refer
the reader to \cite{hirze} for more details, including a description of the one--parameter 
family of generalized Hirzebruch surfaces corresponding to the fundamental triple 
$$(T_a,Q_a,\{(1,0), \pm (0,1),(-1,a) \}).$$

\section{A dictionary}
In recent years, there have been a number of articles on nonrational toric geometry from different viewpoints and new results have been obtained. A common factor of all of these works is of course the presentation of starting convex geometric data. 
A shared feature of these different approaches is the datum of the fundamental triple. Sometimes, further data are added that are instrumental for the constructions; we have already seen an instance of this in Section~\ref{triangolazione}, where the additional data are the ghost vectors. 
Some viewpoints also consider a variant of the notion of quasitorus.

Toric quasifolds can be thought of as stacks, and some authors study nonrational toric geometry within this framework. The first are Hoffman--Sjaamar in the symplectic category \cite[Examples~7.4.3, 9.2.5 and Remark~7.4.4]{HS}.
In their work,  a quasilattice is a crossed Lie module $\partial\colon Q\rightarrow V$, where $Q$ is a finitely generated abelian group and $V$ is a real vector space spanned by $\partial(Q)$. Therefore, the quasilattice $Q$ is no longer a subgroup of $V$, but is surjectively mapped onto one, namely $\partial(Q)$. In this way, there are infinitely many $\partial\colon Q\rightarrow V$ that project onto the same pair $(\partial(Q),V)$.
Notice that the quasilattice $Q$ contains more information than $\partial(Q)$; this is, in this setting, the extra datum that we were referring to earlier. 
Remark also that the pair $(\partial(Q),V)$ gives rise to the quasitorus $V/\partial(Q)$. Here the group that generalizes the torus is called {\em stacky torus}. The definition is quite involved; for our purposes it is sufficient to recall that the datum of a stacky torus $G$ is equivalent to their notion of quasilattice. In \cite[Section~5]{H} Hoffman introduces the starting triple  $(\Delta,G,\Lambda_{f\in F^{max}})$, where $\Delta$ is a simple convex polytope, $G$ is a stacky torus, $F^{max}$ is the set of facets of $\Delta$ and ${\Lambda_f }$ is a free subgroup of rank $1$, given by the intersection of $\partial(Q)$ with the straight line normal to the facet $f$. Notice that there is a unique inward pointing vector that generates ${\Lambda_f }$. This corresponds to the vector, in the fundamental triple, that is associated with the facet $f$. 
We can draw the following diagram 
$$\xymatrix{(\Delta,Q,\{X_1,\ldots,X_d\})\ar@{-->}[r]&(\Delta,G,\Lambda_{f\in F^{max}}),\ar@<1ex>[l]}$$
where ${\Lambda_f}$ is generated by the normal $X_i$ in the corresponding ray and $G$ is any 
$\partial\colon\widetilde{Q} \rightarrow V$ 
with $\partial(\widetilde{Q})=Q$. The arrows have the same meaning as in diagram~(\ref{diagramma}): given a fundamental triple, there are infinitely many triples of the kind 
$(\Delta,G,\Lambda_{f\in F^{max}})$ that project onto it. Viceversa, any triple $(\Delta,G,\Lambda_{f\in F^{max}})$ uniquely determines a fundamental triple. Now, let us recall the following observation in \cite{H}: assigning the triple $(\Delta,G,\Lambda_{f\in F^{max}})$ is equivalent to assigning the triple $(\Delta,G,\Lambda_{f\in F})$, where $F$ is the set of {\em all} faces of $\Delta$ and the groups $\Lambda_f$ are subgroups of $\partial(\widetilde{Q})$ satisfying certain conditions. The triple $(\Delta,G,\Lambda_{f\in F})$ is said to be a {\em decorated stacky moment polytope}, for which a natural notion of isomorphism is given.
The author defines {\it per se} symplectic toric stacks and then proves that the moment mapping defines a bijective correspondence between the set of isomorphism classes of decorated stacky moment polytopes and the set of equivalence classes of symplectic toric stacks \cite[Theorem~6.1]{H}. In our understanding, each equivalence class of symplectic toric stacks corresponds to the symplectic toric quasifold constructed from the associated fundamental triple.

Another article that addresses the problem of generalized toric manifolds is that by Ishida--Krutowsky--Panov \cite{IKP}. 
Their focus is on cohomology. In \cite[Definition~5.5]{IKP} they introduce the quadruple $(\widetilde{V},\widetilde{\Gamma},\widetilde{\Sigma},\widetilde{\lambda})$, 
which they call {\em marked fan}, where $\widetilde{V}$ is a finite--dimensional real vector space, 
$\widetilde{\Gamma}$ is a quasilattice in $\widetilde{V}$, $\widetilde{\Sigma}$ is a fan that is quasirational with respect to $\widetilde{\Gamma}$ and is additionally assumed to be complete and simplicial. Finally, 
$\widetilde{\lambda}$ is a function on the set of one--dimensional cones of $\widetilde{\Sigma}$, with values in the quasilattice 
$\widetilde{\Gamma}$, such that $\widetilde{\lambda}(\rho)$ is a generator of $\rho$. Therefore, $\widetilde{\lambda}$
corresponds exactly to a set of ray generators in the fundamental triple.
The authors then consider compact, connected, complex manifolds with maximal torus actions. These were defined by Ishida in \cite{I2} and later endowed with a canonical foliation in \cite{I1}. In this class of manifolds, the authors define an equivalence relation called  principal equivalence. Then they prove,
building on \cite{I2}, that the set of equivalence classes is in bijective correspondence with isomorphism classes of marked fans \cite[Theorem~5.7]{IKP}. Using this and moment angle manifolds, they are able to drop the hypothesis of shellability in the result by Battaglia--Zaffran \cite{bz1} on the basic cohomology ring of complete simplicial shellable fans. Further results on the cohomology of complete simplicial fans can be found in the recent paper by Krutowsky--Panov \cite{KP}. For the connection between moment angle manifolds and complex manifolds with maximal torus action see also \cite{U}.
As recalled in Section~\ref{triangolazione}, by the construction in \cite{bz1}, we are able to associate with a given fundamental triple a family of LVMB manifolds. Each of these is endowed with a maximal torus action and, therefore, belongs to the equivalence class, defined in \cite{IKP}, corresponding to the given fundamental triple. We know that the leaf space of each of these LVMB manifolds is isomorphic to the complex toric quasifold corresponding to the fundamental triple. We expect that the leaf space of each complex manifold with maximal torus action in the equivalence class is isomorphic to that complex toric quasifold. 

The problem of generalized toric manifolds is again addressed in 
the framework of stacks in a recent article by 
Lupercio--Meersseman--Verjovsky--Katzarkov \cite{KLMV}. 
As starting convex data, they consider {\em quantum fans} and {\em calibrated quantum fans} \cite[Section~4]{KLMV}, which  correspond to fundamental triples and augmented triples, respectively. More specifically, a quantum fan in a quasilattice $\Gamma\subset\mathbb{R}^n$ (which they have renamed quantum lattice) is a pair $(\Delta,v)$, where $\Delta$ is a fan quasirational with respect to $\Gamma$ and $v$ is a set of rays generators contained in $\Gamma$, one for each $1$--dimensional cone of $\Delta$. 
On the other hand, a calibrated quantum fan is a quantum fan plus the additional datum of a calibration. 
Denote by $\{e_1,\ldots, e_p\}$ the standard basis of $\mathbb{Z}^p$. Then a calibration 
is an epimorphism $h\colon \mathbb{Z}^p\longrightarrow \Gamma$ together with a subset $J$ of $\{1,\ldots,p\}$, with the following property: the vectors $h(e_i)$, with $i\notin J$, generate $\mathbb{R}^n$ and give the set $v$ of $d$ rays generators, while the vectors $h(e_i)$, with $i\in J$, are the ghost vectors described in Section~\ref*{triangolazione}. 
The complex quasitorus $\mathbb{C}^n/\Gamma$ is viewed as a quotient stack and called quantum torus.
In correspondence to a complete (calibrated) simplicial quantum fan and to the relative quasilattice $\Gamma$, the authors construct a (calibrated) quantum toric variety. The (calibrated) quantum toric variety is built by suitably gluing the (calibrated) quantum affine toric stacks associated with the maximal cones. From the viewpoint of complex toric quasifolds, this is the affine atlas introduced in \cite[Theorem~2.2]{cx}. They define a notion of morphism for (calibrated) quantum fans and (calibrated) quantum toric varieties. The map between these two categories given by the construction turns out to be functorial, it is naturally surjective, and is proved to be an equivalence of categories 
\cite[Theorems~5.18,~6.24]{KLMV}.

Finally, Boivin extends the above equivalence of categories to calibrated nonsimplicial fans and the corresponding calibrated quantum toric varieties \cite[Theorem~4.2.2.2]{boivin}. He makes use of the same starting convex data, that he calls (calibrated) quantum fans and quantum lattices as well. But he needs to introduce an auxiliary datum, a further calibration, in order
to deal with nonsimpliciality.

The following diagram gives a synthetic and unified picture of the various constructions:
$$	
\xymatrix@M=6pt@R=50pt@C=34pt{*+[F-]{\txt<12pc>{triple + additional data}}\ar[r]^{\quad\overset{\approx}{\Psi}}\ar[d]&
*+[F-]\txt<9pc>{geometric space}\ar[d]\\
\txt{fundamental triple}\ar[r]^{\quad\Psi}\ar[ur]^{\quad\quad\widetilde{\Psi}}&\txt{toric quasifold}
}
$$
The box on the left represent a whole family of data projecting down to the same fundamental triple, the box on the right represent the corresponding family of geometric spaces.
The lower level is the mapping that associates the fundamental triple with the toric quasifold.
Notice that we can always construct, from a given fundamental triple, an object that lies in the family above.
The articles \cite{bz1,H} and, in the calibrated case, \cite{KLMV,boivin} can be viewed as instances of the mapping $\overset{\approx}{\Psi}$, \cite{IKP} can be viewed as an instance of the mapping $\widetilde{\Psi}$, and finally \cite{pcras,p,cx}, and \cite{KLMV,boivin} can be viewed as instances of the mapping $\Psi$.
In our understanding, in each of the above--mentioned constructions relative to the mappings $\overset{\approx}{\Psi}$ and $\widetilde{\Psi}$, there is a mapping that projects any space of the family to the toric quasifold. An instance is given in \cite[Theorem~2.1]{rome}.

We conclude by briefly mentioning a number of other related works. Ratiu--Zung \cite{ratiu} and Lin--Sjamaar \cite{LS} study nonrational convex polytopes in the context of presymplectic manifolds. The first authors specialize to the toric case, but a triple is not explicitly provided. Pir--Sottile \cite{PS} introduce the notion of {\em irrational toric variety} for arbitrary fans and show that, when the fan is the normal fan of a polytope, the irrational toric variety is homeomorphic to that polytope. Quasifolds have been studied in the framework of diffeology in \cite{IZP}. We expect that this approach will have implications in nonrational toric geometry. Finally, we remark that Bressler--Lunts \cite{BL1, BL2} and Karu \cite{karu} devised a powerful approach to the combinatorics of nonrational polytopes that extends cohomological properties of toric varieties to general fans, without constructing a corresponding toric space.

\bigskip

{\small 

\noindent
Dipartimento di Matematica e Informatica "U. Dini", Universit\`a di Firenze \\
Viale Morgagni 67/A, 50134 Firenze, ITALY

\noindent
fiammetta.battaglia@unifi.it, elisa.prato@unifi.it}

\end{document}